\documentclass[10pt]{amsart}

\usepackage{amssymb, amscd, amsmath, amsthm,  graphicx, latexsym}

\def\zz{{\bf Z}}

\def\C{{\bf C}}

\newtheorem{theorem}{Theorem}

\newtheorem{corollary}[theorem]{Corollary}

\begin{document}

\title[Ozsv\'{a}th--Szab\'{o} and   Rasmussen   invariants of cable knots]
{Ozsv\'{a}th--Szab\'{o} and   Rasmussen   invariants of cable knots}

\author{Cornelia A. Van Cott}
\address{
University of San Francisco, San Francisco, California 94117} 
\email{cvancott@usfca.edu}
\maketitle

\begin{abstract} 
We study the behavior of the Ozsv\'{a}th--Szab\'{o} and Rasmussen knot concordance invariants $\tau$ and $s$ on $K_{m,n}$, the $(m,n)$-cable of a knot $K$ where $m$ and $n$ are relatively prime.  We show that for every knot $K$ and for any fixed positive integer $m$, both of the invariants evaluated on $K_{m,n}$ differ from their value on the torus knot $T_{m,n}$ by fixed constants for all but finitely many $n>0$.  Combining this result together with Hedden's extensive work on the behavior of $\tau$ on $(m,mr+1)$-cables yields bounds on the value of $\tau$ on any $(m,n)$-cable of $K$.  In addition, several of Hedden's obstructions for cables bounding complex curves are extended.
\end{abstract}

\section{Introduction}  
The $(m,n)$-cable of a knot $K$, denoted $K_{m,n}$, is the satellite knot with companion $K$ and pattern $T_{m,n}$, the $(m,n)$-torus knot.  The behavior of many classical concordance invariants has been shown to be rather predictable on cable knots.  For example, it is a classical result (c.f.~\cite{lic}) that the Alexander polynomial of a cable knot is given by 
$$\Delta_{K_{m,n}}(t) = \Delta_K(t^m)\Delta_{T_{m,n}}(t).$$  
Shinohara~\cite{shi} found a formula for the signature of a cable knot, and Litherland~\cite{lith} extended the result, finding the value of Tristam-Levine signatures on a cable knot to be
$$\sigma_\omega(K_{m,n}) = \sigma_{\omega^m}(K) + \sigma_\omega(T_{m,n}).$$
Milnor signatures and Casson-Gordon invariants of cables (see~\cite{lithcg} and~\cite{kear}, respectively, for details) also yield nice formulas.  

The purpose of this note is to investigate two relatively new concordance invariants --  the Ozsv\'{a}th--Szab\'{o} invariant $\tau$ and the Rasmussen invariant $s$  -- and their behavior on cable knots.  The discussion here will use only the formal properties that the two invariants have in common.   

Both $\tau$ and $s$ were introduced in connection with developments in the theory of knot homologies: $\tau$ is defined in terms of knot Floer homology~\cite{os, ra} and the Rasmussen invariant $s$ is defined in terms of Khovanov homology~\cite{ra}.    These two invariants have enabled important progress in the field of knot theory, providing new proofs for Milnor's conjecture~\cite{os, ra} and examples of Alexander polynomial one knots which are not smoothly slice~\cite{liv}.


No work has been done to compute the Rasmussen invariant for cables, but the behavior of the Ozsv\'{a}th--Szab\'{o} concordance invariant $\tau$ under $(m, mr+1)$-cabling has been investigated by Hedden~\cite{h1,h2}.  Through careful investigation of the relationship between the filtered chain homotopy types of $\mathcal{F}(K_{m,mr+1},i)$ and $\mathcal{F}(K,i)$, he obtained the following main result:\\

\begin{theorem}~\cite{h2}~\label{matt}
Let $K\subset S^3$ be a nontrivial knot.  Then the following inequality holds for all $r$:
$$m\tau(K) + \frac{(mr)(m-1)}{2} \leq \tau(K_{m,mr+1}) \leq m  \tau(K) + \frac{(mr)(m-1)}{2} + m-1.$$
In the special case when $K$ satisfies $\tau(K) = g(K)$, we have the equality,
$$\tau(K_{m,mr+1}) = m\tau(K) + \frac{(mr)(m-1)}{2},$$
whereas when $\tau(K) = -g(K)$, we have
$$\tau(K_{m,mr+1}) = m\tau(K) + \frac{(mr)(m-1)}{2} + m-1.$$
\end{theorem}

When appropriately normalized, $\tau$ and $s$ share several formal properties and agree on many families of knots, though in general they have been shown to be distinct invariants~\cite{ho}.  Stated in reference to $\tau$, the essential formal properties are as follows~\cite{os2, ra}:
\begin{enumerate}
\item $\tau$ is a homomorphism from the smooth knot concordance group $\mathcal{C}$ to $\zz$. 
\item $|\tau(K)| \leq g_4(K)$, where $g_4(K)$ denotes the 4-genus of $K$.
\item $\tau(T_{m,n}) = \frac{(m-1)(n-1)}{2},$ where $T_{m,n}$ denotes the $(m,n)$-torus knot with $m,n \geq 1$.
\end{enumerate}

It can be shown that $s/2$ also satisfies these three properties~\cite{ra2}.  Our main results will only depend on these formal properties, and hence apply to both invariants.  To proceed concisely, let $\nu$ denote any concordance invariant satisfying the above three properties.  

Fixing $m>0$, we would like to study the value of $\nu$ on $K_{m,n}$ as a function of $n$, where $n$ ranges over the integers relatively prime to $m$.  (Notice that $K_{m,n} = - K_{-m,-n}$, and so the restriction $m > 0$ does not limit our results.)  From our observations about other concordance invariants, we expect that the behavior of $\nu(K_{m,n})$ as a function of $n$ is somehow related to the behavior of $\nu(T_{m,n})$.  This, in fact, is true.  As a function of $n$, $\nu(T_{m,n})$ is linear of slope $\frac{(m-1)}{2}$ for $n>0$.  We will see that the function $\nu(K_{m,n})$ is close to being linear with the same slope.  Specifically, we subtract from $\nu$ a linear function to construct the following function:
$$h(n) = \nu(K_{m,n}) - \frac{(m-1)}{2}n,$$
where $n$ is an integer relatively prime to $m$.  We have the following theorem:

\begin{theorem}\label{me}
The function $h(n)$ is a nonincreasing $\frac{1}{2}\cdot \zz$-valued function which is bounded below. In particular, we have
$$-(m-1) \leq h(n) - h(r) \leq 0$$
for all $n > r$, where both $n$ and $r$ are relatively prime to $m$.
\end{theorem}

From this result it follows that for all $n$ large enough, $h$ is constant.  Hence for $n$ large enough, $\nu(K_{m,n})$ differs from $\nu(T_{m,n})$ by a fixed constant.  That is, for every knot $K$ there exist integers $N$ and $c$ such that $\nu(K_{m,n}) = \nu(T_{m,n}) + c$ for all $n > N$, where $n$ is relatively prime to $m$.  Additionally, a similar statement with corresponding constant $c'$ holds for all $n< N'$ for some $N'$.

Theorem~\ref{me} is sharp in the sense that there are knots $K$ with associated functions $h$ which achieve the bounds given in the theorem.  For example, when $K$ is slice,  $h(n) = (m-1)/2$ for all $n < 0$ and $h(n) = -(m-1)/2$ for all $n > 0$.  Here the drop in functional value from $n = -1$ to $n = 1$ is maximal:  $h(1) - h(-1) = -(m-1)$.  On the other hand, we will see that when $\nu = \tau$ and $\tau(K) = g_3(K)$, the function $h$ is constant.  

Using Theorem~\ref{me}, we can take several results which apply only to $(m,mr+1)$-cables and extend their scope to include {\em all} cables.  For example, the bounds on the value of $\tau$ on $(m,mr+1)$-cables described in Theorem~\ref{matt} extend to all cables as follows.
\begin{corollary}~\label{fill}
Let $K\subset S^3$ be a nontrivial knot.  Then the following inequality holds for all $n$ relatively prime to $m$:
$$m\tau(K) + \frac{(m-1)(n-1)}{2} \leq \tau(K_{m,n}) \leq m  \tau(K) + \frac{(m-1)(n+1)}{2}.$$
When $K$ satisfies $\tau(K) = g(K)$, we have $\tau(K_{m,n}) = m\tau(K) + \frac{(m-1)(n-1)}{2},$
whereas when $\tau(K) = -g(K)$, we have $\tau(K_{m,n}) = m\tau(K) + \frac{(m-1)(n+1)}{2}.$
\end{corollary}
 
Observe that the results in Corollary~\ref{fill} could probably also have been obtained by using the definition of $\tau$ and studying the filtered chain homotopy  type of $\mathcal{F}(K_{m,n})$ for $n$ relatively prime to $m$.  However, the proof here avoids this and uses only the analysis of $\mathcal{F}(K_{m,mr+1})$ in~\cite{h2} together with Theorem~\ref{me} to obtain the result for all cables.

The second half of Corollary~\ref{fill} motivates studying knots $K$ for which $\tau(K) = g(K)$.  Hedden summarized many results about such knots and their $(m,mr+1)$-cables in~\cite{h2}.  Now combining that discussion with Corollary~\ref{fill} from above, we can extend several of his results to a more general setting.  Let $\mathcal{P}$ denote the class of all knots satisfying the equality $\tau(K) = g(K)$.  An immediate consequence of Corollary~\ref{fill} is the following.

\begin{corollary}
Let $K$ be a nontrivial knot in $S^3$, and let $n$ be relatively prime to $m$.
\begin{enumerate}
\item If $K \in \mathcal{P}$, then $K_{m,n}\in \mathcal{P}$ if and only if $n>0$. 
\item If $K \notin \mathcal{P}$, then $K_{m,n}\notin \mathcal{P}$.
\end{enumerate}
\end{corollary}

As discussed in~\cite{h2}, $\mathcal{P}$ contains several classes of knots.  We mention two such classes here:
\begin{itemize}
\item Any knot $K$ which bounds a properly embedded complex curve, $V_f \subset B^4$, with $g(V_f) = g(K)$. This set of knots includes, for example, positive knots (that is, knots which admit diagrams with only positive crossings).~\cite{h0,liv}
\item Any knot which admits a positive lens space (or L-space) surgery.~\cite{os3}
\end{itemize}
From this, we have the following immediate applications extending the work of~\cite{h2}.

\begin{corollary}
If $K_{m,n}$ bounds a properly embedded complex curve $V_f \subset B^4$ satisfying $g(V_f) = g(K_{m,n})$, then $n >0$ and $\tau(K) = g(K)$.
\end{corollary}

\begin{corollary}
Suppose that $K_{m,n}$ admits a positive lens space (or L-space) surgery.  Then $n > 0$ and $\tau(K) = g(K)$.
\end{corollary}
 
 \begin{corollary}
 Suppose $K \notin \mathcal{P}$.  Then $K_{m,n}$ is not a positive knot for any relatively prime pair of integers $m$, $n$.
 \end{corollary}
 
 A final corollary concerns a more general class of knots --- the class of $\C$--knots.  A knot $K$ is a $\C$--knot if $K$ bounds a properly embedded complex curve $V_f \subset B^4$.  From~\cite{h0,pl,rud}, we know that for such knots, $\tau(K) = g_4(K) \geq 0$.  Coupling this result with Corollary~\ref{fill}, we have the following corollary.
\begin{corollary}
Suppose that $K_{m,n}$ is a $\C$-knot.  Then $n \geq \frac{-2m\tau(K)}{m-1}-1$.
\end{corollary}
 
The primary significance of each of these corollaries is that they can be used as obstructions to cables having the discussed properties.  Moreover, it is interesting that $\tau$ provides obstructions to such a wide array of  geometric notions.  For an excellent extended discussion of this, we refer the reader to~\cite{h2}.
 
This paper is organized as follows.  Section~\ref{mine} contains the proof of Theorem~\ref{me}.  Section~\ref{apps} contains the proof of Corollary~\ref{fill}.  Finally, in Section~\ref{braids} we observe that the strategy for the proof of Theorem~\ref{me} extends to a broader setting in which, instead of cabling, we consider a braiding construction.  

\subsection*{Acknowledgments} I thank both Chuck Livingston and Matt Hedden for several helpful conversations.  

\section{Proof of Theorem~\ref{me}}~\label{mine}

Let $r, n$ be integers relatively prime to $m$ with $n > r$.  The general strategy here is to first find a cobordism between $K_{m,n} \# -K_{m,r}$ and a torus knot.  
 
We begin with the knot $K_{m,n} \# -K_{m,r}$.  Working through signs and orientations carefully, we find that 
$$K_{m,n} \# -K_{m,r} = K_{m,n} \# (-K)_{m,-r}.$$  

We will now do a series of band moves to the knot $K_{m,n} \# (-K)_{m,-r}$.  A band move on any knot $K \subset S^3$ is accomplished as follows.  Start with an embedding $b:I\times I \longrightarrow S^3$ such that $b(I\times I)\cap K = b(I\times \{0,1\})$ and such that $b$ respects the orientation of $K$.  Define $K_b = K - b(I\times\{0,1\}) \cup b(\{0,1\} \times I)$.  The knot (or link) $K_b$ is the result of doing a band move along $b$.  Doing a band move to a knot simultaneously constructs a cobordism from the knot $K$ to $K_b$. The genus of this cobordism can be computed explicitly.  For example, in the special case that the result of performing a sequence of band moves is again a knot, one can show that the genus of the cobordism is: $\frac{1}{2} \cdot $number of bands added.

\begin{figure}
\begin{center}
\includegraphics[width=10cm]{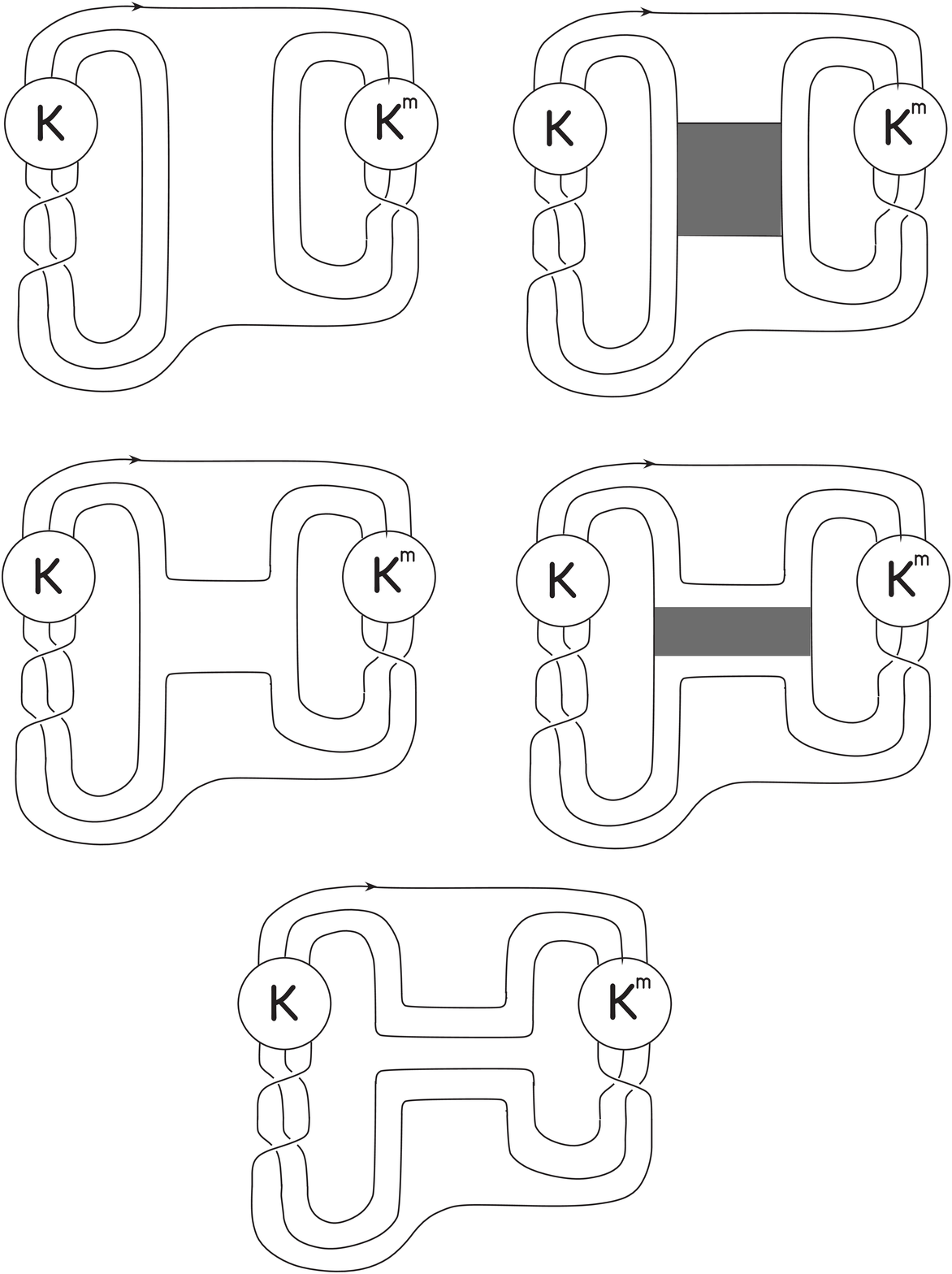}
\end{center}
\caption{Beginning with the knot $K_{3,2} \# (-K)_{3,-1}$, we perform two band moves and obtain the knot $(K \# -K)_{3, 1}$.  $K^m$ denotes the mirror image of $K$.}\label{bandmoves}
\end{figure}

\begin{figure}
\begin{center}
\includegraphics[width=10cm]{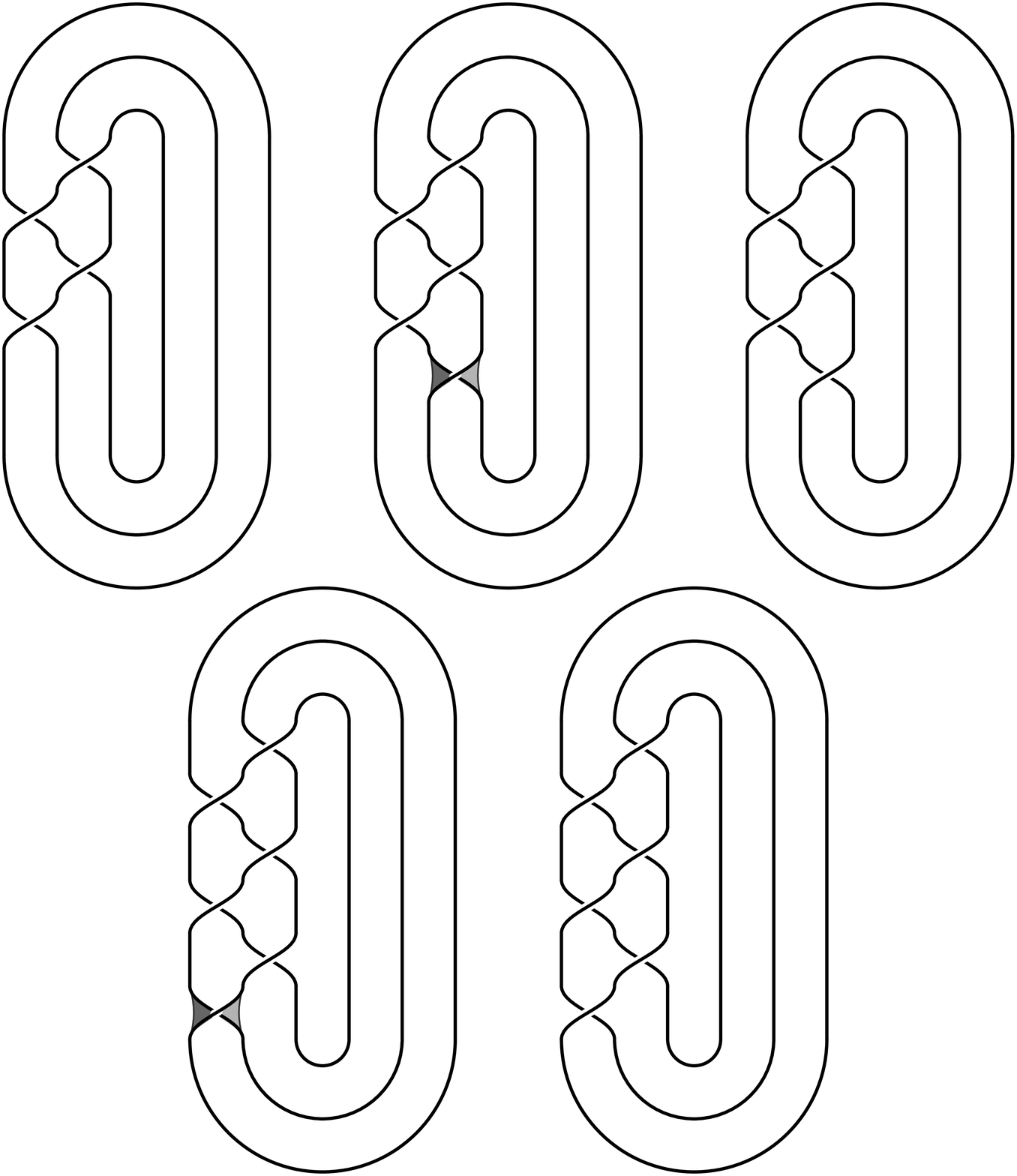}
\end{center}
\caption{Beginning with the torus knot $T_{3,2}$, we perform two band moves and obtain $T_{3,3}$.}\label{torusbands}
\end{figure}

Now there is a sequence of $m-1$ band moves on $K_{m,n} \# (-K)_{m,-r}$ which results in the knot (or link) $(K \# -K)_{m, n-r}$.  See Figure~\ref{bandmoves} for an example.  Since $K\#-K$ is cobordant to the unknot, $(K \# -K)_{m, n-r}$ is cobordant to the torus link $T_{m,n-r}$.  Let $k_+$ denote the smallest positive integer such that $n-r+k_+$ is relatively prime to $m$.  (If $n-r$ is already relatively prime to $m$, then set $k_+ = 0$.)  By doing $k_+\!\cdot\!(m-1)$ additional band moves to the torus link $T_{m,n-r}$, we obtain the torus knot $T_{m, n-r+k_+}$ (Figure~\ref{torusbands}).  Altogether, the total number of band moves performed was $(k_++1)(m-1)$.  Therefore, the knot $K_{m,n} \# -K_{m,r}$ is genus $(k_++1)(m-1)/2$ cobordant to the torus knot $T_{m, n-r+k_+}$.  Hence we conclude that 
$$g_4(K_{m,n} \# -K_{m,r} \# -T_{m, n-r+k_+}) \leq (k_++1)(m-1)/2.$$
Now since $|\nu(K)|\leq g_4(K)$, it follows that
$$|\nu(K_{m,n} \# -K_{m,r} \# -T_{m, n-r+k_+})| \leq (k_++1)(m-1)/2.$$
Simplifying the expression using the properties of $\nu$, we obtain

$$\Big|\nu(K_{m,n}) -\nu(K_{m,r})-\frac{(m-1)(n-r+k_+-1)}{2}\Big| \leq \frac{(k_++1)(m-1)}{2}.$$

At this point, recall the function $h(n)$ which we defined earlier.  Using the definition of $h$, we can further simplify the inequality:
$$\Big|h(n) - h(r) - \frac{(m-1)(k_+-1)}{2}\Big| \leq \frac{(k_++1)(m-1)}{2}.$$
Hence, 
\begin{equation}\label{1}
-(m-1) \leq h(n) - h(r) \leq k_+(m-1).
\end{equation}
Notice that if $k_+=0$, then we are done.  If not, then we continue as follows.

Similar to before, let $k_-$ denote the largest negative integer such that $n-r+k_-$ is relatively prime to $m$.  By doing $|k_-|\!\cdot\!(m-1)$ band moves to $T_{m,n-r}$, we can obtain the torus knot $T_{m, n-r+k_-}$.  Proceeding through the same steps as before, we obtain
 \begin{equation}\label{2}
(k_--1)(m-1) \leq h(n) - h(r) \leq 0.
 \end{equation}
 
Combining (\ref{1}) and (\ref{2}), we have 
$$-(m-1) \leq h(n) - h(r) \leq 0$$
for all integers $n > r $ where both $n$ and $r$ are relatively prime to $m$. 

\section{Proof of Corollary~\ref{fill}}~\label{apps}
Combining Theorem~\ref{matt} and Theorem~\ref{me} together, we obtain an easy proof that the bounds on the value of $\tau$ on $(m,mr+1)$-cables described in Theorem~\ref{matt} extend to all cables.  We now restate and prove Corollary~\ref{fill}.
\\

\noindent{\bf Corollary 3.}
{\em Let $K\subset S^3$ be a nontrivial knot.  Then the following inequality holds for all $n$ relatively prime to $m$:
$$m\tau(K) + \frac{(m-1)(n-1)}{2} \leq \tau(K_{m,n}) \leq m  \tau(K) + \frac{(m-1)(n+1)}{2}.$$
When $K$ satisfies $\tau(K) = g(K)$, we have 
$$\tau(K_{m,n}) = m\tau(K) + \frac{(m-1)(n-1)}{2},$$
whereas when $\tau(K) = -g(K)$, we have 
$$\tau(K_{m,n}) = m\tau(K) + \frac{(m-1)(n+1)}{2}.$$
}

\begin{proof}
The proof of this corollary is obtained by carefully combining the equalities and inequalities found in Theorem~\ref{matt} and Theorem~\ref{me}.  We will demonstrate a portion of the proof, leaving the rest to the reader.  

Let $m$ and $n$ be two relatively prime integers with $m>0$.  Let $r$ be an integer such that $n > mr+1$.  Then by Theorem~\ref{me}, 
$$h(n) - h(mr+1) \leq 0.$$
Using the definition of $h$ and letting $\nu = \tau$, we obtain
$$\tau(K_{m,n}) \leq \tau(K_{m,mr+1}) - \frac{m-1}{2}(mr - n + 1).$$
Using the upper bound on $\tau(K_{m,mr+1})$ given by Theorem~\ref{matt}, we have
$$\tau(K_{m,n}) \leq m\tau(K) + \frac{(m-1)(n+1)}{2},$$
which is one side of the desired inequality.  

To obtain the other side of the inequality, let $r'$ be an integer such that $mr'+1 > n$.  Then by Theorem~\ref{me}, 
$$h(mr'+1) - h(n) \leq 0.$$
We leave to the reader the task of reducing this inequality (using methods exactly similar to above) to obtain the desired second half of the inequality in the corollary.

Now let $K$ be a knot such that $\tau(K) = g(K)$.  Suppose for contradiction that  $\tau(K_{m,n}) \neq m\tau(K) + \frac{(m-1)(n-1)}{2}.$  By the inequality discussed above, this implies that  $\tau(K_{m,n}) > m\tau(K) + \frac{(m-1)(n-1)}{2}.$  Again, let $r$ be an integer such that $n > mr + 1$.  Then we have
\begin{eqnarray*}
h(n) - h(mr+1)  &=&  \tau(K_{m,n}) - \frac{(m-1)}{2}n - \tau(K_{m,mr+1}) + \frac{m-1}{2}(mr+1)\\
&=&  \tau(K_{m,n}) - \frac{(m-1)}{2}n - m\tau(K) + \frac{(m-1)}{2}\\
&>& m\tau(K) + \frac{(m-1)(n-1)}{2} - \frac{(m-1)}{2}n - m\tau(K) + \frac{(m-1)}{2}\\
&=& 0.
\end{eqnarray*}

This is contradicts Theorem~\ref{me}.  Therefore, $\tau(K_{m,n}) = m\tau(K) + \frac{(m-1)(n-1)}{2}$ for all $n$ relatively prime to $m$.  A similar argument settles the case when $K$ is a knot such that $\tau(K) = -g(K)$.
\end{proof}

\section{Further analysis}~\label{braids}
The process of cabling a knot can be reinterpreted as a special case of the following more general procedure.   Let $\beta$ be an element of the braid group $B_m$ such that the closure of the braid $\widehat \beta$ is a knot.  There is a natural solid torus $V$ which contains the closed braid $\widehat \beta$.  Remove a neighborhood of a knot $K$ in $S^3$ and glue in the solid torus $V$ by a homeomorphism which maps longitude to longitude and meridian to meridian.  We denote the resulting knot by $K_\beta$.  Notice that if we take the braid $\beta \in B_m$ to be $(\sigma_{m-1}\sigma_{m-2}\cdots\sigma_1)^n$ (where $\sigma_i$ denotes the $i^{th}$ standard generator of the braid group), then the resulting knot $K_\beta$ is the $(m,n)$-cable $K_{m,n}$.  

For any braid $\beta \in B_m$, let $\beta_r$ denote the braid consisting of $\beta$ with $r$ full twists adjoined to the end of the braid.  Specifically, $\beta_r = \beta(\sigma_{m-1}\sigma_{m-2}\cdots\sigma_1)^{mr}$.  The value of $\nu$ on $K_{\beta_r}$ as a function of $r$ turns out to have controlled behavior similar to that of cabling.  Define the function 
$$g(r) = \nu(K_{\beta_r}) - \frac{(m-1)}{2}mr,$$
where $\beta \in B_m$ is a braid whose closure is a knot and $r$ is an integer.  Then we have the following theorem about the behavior of the function $g$.
\\
\begin{theorem}~\label{braiding}
The function $g(r)$ is a nonincreasing integer valued function which is bounded below. In particular, 
$$-(m-1) \leq g(r) - g(s) \leq 0$$
for all $r > s$.
\end{theorem}

From this theorem, it follows that the function $g$ is eventually constant.  This allows us to describe quite clearly a relationship among the values of $\tau$ (and $s$) on an entirely new set of knots.  Fixing a knot $K$ and a braid $\beta \in B_m$ such that $\widehat{\beta}$ is a knot, Theorem~\ref{braiding} implies that for all large $r$, 
$$\nu(K_{{\beta}_{r+1}}) = \nu(K_{{\beta}_{r}}) + \frac{m(m-1)}{2}.$$
where $\nu$ can be taken to be either $\tau$ or $s$. 
Note that if we take $K$ in the above construction to be the unknot, then the theorem relates the values of $\nu$ on knots with braid representatives which differ by full twists.  

We turn now to the proof of Theorem~\ref{braiding}.

\begin{proof}
As with the proof of Theorem~\ref{me}, the first goal here is to find a cobordism between $K_{\beta_r} \# -K_{\beta_s}$ and a torus knot.  Notice that $-K_{\beta_s} = (-K)_{(\beta^{-1})_{-s}}$.  Therefore, 
$$K_{\beta_r} \# -K_{\beta_s} = K_{\beta_r} \# (-K)_{(\beta^{-1})_{-s}}.$$   
By doing $m-1$ band moves to the latter knot, we obtain the knot $(K \# -K)_{(\beta\beta^{-1})_{r-s}}$.  Since $K\#-K$ is cobordant to the unknot and $\beta\beta^{-1}$ is the trivial $m$-strand braid, this new knot is cobordant to the torus link $T_{m,m(r-s)}$.  Again, by doing $(m-1)$ band moves to the torus link $T_{m,m(r-s)}$, we obtain the torus knot $T_{m, m(r-s)+1}$.  A total of $2(m-1)$ band moves have been performed.  Therefore, the knot $K_{\beta_r} \# -K_{\beta_s}$ is genus $(m-1)$ cobordant to the torus knot $T_{m, m(r-s)+1}$.  Hence 
$$g_4(K_{\beta_r} \# -K_{\beta_s} \# -T_{m, m(r-s)+1}) \leq m-1.$$
Now since $|\nu(K)|\leq g_4(K)$, it follows that
$$|\nu(K_{\beta_r} \# -K_{\beta_s} \# -T_{m, m(r-s)+1})| \leq m-1,$$
which simplifies to

$$\Big|\nu(K_{\beta_r}) -\nu(K_{\beta_s})-\frac{(m-1)m(r-s)}{2}\Big| \leq m-1.$$

We now recall the function $g(r)$ which we defined earlier.  Using the definition of $g$, we can further simplify the inequality and obtain:
\begin{equation}\label{3}
-(m-1) \leq g(r) - g(s) \leq m-1.
\end{equation}
This gives us only half of the desired inequality.  To obtain the remaining half, go back to the torus link $T_{m,m(r-s)}$ which we obtained from $K_{\beta_r} \# -K_{\beta_s}$ by a cobordism which added $m-1$ bands.  Instead of adding $m-1$ additional bands to obtain the torus knot $T_{m,m(r-s)+1}$, add  $m-1$ bands to obtain the torus knot $T_{m,m(r-s)-1}$.    Proceeding through the same steps as before, we obtain
 \begin{equation}\label{4}
-2(m-1) \leq g(r) - g(s) \leq 0.
 \end{equation}
 
Combining (\ref{3}) and (\ref{4}), we have 
$$-(m-1) \leq g(r) - g(s) \leq 0$$
for all integers $r > s$, as desired. 
\end{proof}


\begin{thebibliography}{9999}
 
 
 
\bibitem[1]{h0} {\bf {M.~Hedden}}, {\sl Notions of positivity and the Ozsv\'{a}th--Szab\'{o} concordance invariant,}  arxiv.org/math/0509499.
 
\bibitem[2]{h1} {\bf {M.~Hedden}}, {\sl On knot Floer homology and cabling,}  Algebr. Geom. Topol. {\bf 5} (2005) 1197--1222.
 
\bibitem[3]{h2}{\bf {M.~Hedden}}, {\sl On knot Floer homology and cabling II,}  	arXiv:0806.2172v2.
 
\bibitem[4]{ho} {\bf M.~Hedden}, {\bf P.~Ording}, {\sl The Ozsv\'{a}th--Szab\'{o} and Rasmussen concordance invariants are not equal},  arxiv.org/math/0512348.
 
\bibitem[5]{kear}{\bf C.~Kearton}, {\sl The Milnor signatures of compound knots}, Proc. Amer. Math. Soc. {\bf 76} (1979), no.~1, 157--160. 
 
 
 \bibitem[6]{lic}{\bf W.~B.~R.~ Lickorish}, {\sl An introduction to knot theory}. {\sl Graduate Texts in Mathematics}, {\bf 175}. Springer-Verlag, New York, 1997.
 
\bibitem[7]{liv} {\bf C.~Livingston}, {\sl Computations of the Ozsv\'{a}th--Szab\'{o} knot concordance invariant}, Geom. Topol. {\bf 8} (2004) 735--742.

\bibitem[8]{lith}{\bf R.~A.~Litherland}, Signatures of iterated torus knots. {\sl Topology of low-dimensional manifolds (Proc. Second Sussex Conf., Chelwood Gate, 1977), } Lecture Notes in Math., {\bf 722}, Springer, Berlin, 1979, 71--84. 

\bibitem[9]{lithcg}{\bf R.~A.~Litherland}, Cobordism of satellite knots. {\sl Four-manifold theory} (Durham, NH, 1982), {\sl Contemp. Math.}, {\bf 35}, Amer. Math. Soc., Providence, RI, 1984, 327--362.

 
 \bibitem[10]{os} {\bf P. ~Ozsv\'{a}th}, {\bf  Z. ~Szab\'{o}}, {\it  Heegaard Floer homology and alternating knots},  Geom. Topol.  {\bf  7} (2003)  225--254.    
 
\bibitem[11]{os2} {\bf P. ~Ozsv\'{a}th}, {\bf  Z. ~Szab\'{o}}, {\it  Knot Floer homology and the four-ball
genus},  Geom. Topol.  {\bf  7} (2003)  615--639.    

\bibitem[12]{os3} {\bf P. ~Ozsv\'{a}th}, {\bf  Z. ~Szab\'{o}}, {\it  On knot Floer homology and lens space surgeries},  Topology.  {\bf  44} (2005), no. 6,  1281 -- 1300.    

\bibitem[13]{pl} {\bf O.~Plamenevskaya},  {\sl Bounds for the Thurston--Bennequin number from Floer homology}, Algebr. Geom. Topol. {\bf 4} (2004), 399--406.


\bibitem[14]{ra} {\bf J.~A.~Rasmussen}, {\sl  Floer homology and knot complements,} PhD thesis, Harvard University, 2003.

\bibitem[15]{ra2} {\bf J.~A.~Rasmussen}, {\sl  Khovanov homology and the slice genus,} arxiv.org/math/0402131.

\bibitem[16]{rud} {\bf L.~Rudolph}, {\sl Algebraic functions and closed braids}, Topology, {\bf 22}, (1983), no.2, 191--202. 

\bibitem[17]{shi} {\bf Y.~Shinohara},  {\sl On the signature of knots and links}, Trans. Amer. Math. Soc.  {\bf 156} (1971), 273--285.

 
\end{thebibliography}
\end{document}